\documentclass[12pt]{article}

\usepackage{latexsym,amssymb,color, hyperref}

\def\tr{^{\rm T}}
\def\Real{\mathbb R}
\def\Compl {\mathbb C}
\def\Integer{\mathbb N}

\def\dst{\displaystyle}
\def\ca{{{\cal A}}}
\def\cb{{{\cal B}}}

\def\qmx#1{\left(\matrix{#1}\right)}

\def\beeq#1{\begin{equation}{#1}\end{equation}}

\def\ba{\begin{array}}
\def\ea{\end{array}}
\def\eqa{\begin{eqnarray}}
\def\eqe{\end{eqnarray}}

\newtheorem{theorem}{Theorem}

\newtheorem{proposition}{Proposition}

\newtheorem{lemma}{Lemma}

\newenvironment{proof}{\medskip\noindent{\it Proof. }}{ \medskip}
\newenvironment{remark}{\medskip\noindent{\it Remark. }}{
\medskip}

\newcounter{comment}
\newlength\comwidth

\begin{document}

\title{Output Stabilization \\via Nonlinear Luenberger Observers \thanks{This work was
partially supported by MIUR and Institute of Advanced Study (ISA) of University
of Bologna, by NSF under grant ECS-0314004.}}

\author{
L. Marconi $^{\circ}$, L. Praly $^{\ast}$, A. Isidori $^{\dag
\circ}$}
\date{}
\maketitle
\begin{center}

$^{\circ}$ C.A.SY. -- D.E.I.S., University of Bologna, Italy.
\smallskip


$^{\ast}$ CAS, \'Ecole des Mines de Paris, Fontainebleau, France.
\smallskip

$^{\dag}$ D.I.S, Universit\`{a} di Roma ``La Sapienza'', Italy and\\
E.S.E., Washington University, St. Louis, USA.

\smallskip
\normalsize
\end{center}


\begin{abstract}
The present paper addresses the problem of existence of an (output) feedback
law to the purposes of asymptotically steering to zero a given controlled
variable, while keeping all state variables bounded, for any initial conditions
in a given compact set. The problem can be viewed as an extension of the
classical problem of semi-globally stabilizing the trajectories of a controlled
system to a compact set. The problem also encompasses a version of the
classical problem of output regulation. Assuming only the existence of a
feedback law that keeps the trajectories
  of the zero dynamics of the controlled plant bounded, it is shown
that there exists a controller solving the problem at hand.
  The paper is deliberately focused on theoretical results regarding the existence
  of such controller. Practical aspects
  involving the design and the implementation of the controller are left to a
  forthcoming work.
\end{abstract}


\section{Introduction} \label{sec1}
 The problem of controlling a system in such a way that some
selected variables converge to zero while all other state
variables remain bounded is a relevant problem in control theory.
It includes, as special cases, the problem of asymptotic
stabilization of a fixed equilibrium point and the problem of
asymptotic stabilization of a fixed invariant set. It also
includes design problems in which some selected variables are
required to asymptotically track (or to asymptotically reject)
certain signals generated by an independent autonomous system.
Problems of this kind, usually referred to as problems of ``output
regulation", have been extensively studied in the past for linear
systems (see \cite{Da, FW76, Fr}) as well as, beginning with the
seminal work \cite{IB}, for nonlinear systems. As a matter of
fact, these problems can be viewed as problems in which a
``regulated" output of an ``augmented system" (a system consisting
of  the controlled plant and of the exogenous system generator)
must be asymptotically steered to zero while all other state
variables are kept bounded. As for instance pointed out in
\cite{IB}, the basic challenges in a problem of this type are to
create an invariant set on which the desired regulated output
vanishes, and  to render this set asymptotically attractive.

Even though limited in scope (the design method suggested therein
being only meant to secure local, and non-robust, regulation about
an equilibrium point) paper \cite{IB}  had the merit of
highlighting a few basic concepts and ideas which shaped all
subsequent developments in this area of research. These ideas
include the fundamental link between the problem in question and
the notion of ``zero dynamics" (a concept introduced and studied
earlier by the same authors), the necessity of the existence of a
(controlled) invariant set on which the desired regulated output
vanishes, and an embryo of design philosophy based on the idea of
making this invariant set locally (and exponentially) attractive.

In the past fifteen years, the design philosophy of \cite{IB} was
extended in several directions. One clear need was to move from
``local" to ``non-local" convergence, a goal which was pursued --
for instance --  in  \cite{Kh2}, \cite{IA} and \cite{SIM99a},
\cite{BI03}, where different approaches (at increasing levels of
generality) have been proposed. Another concern was to obtain
design methods which are insensitive, or even robust, with respect
to model uncertainties. This issue was originally addressed in
\cite{HL}, where it was shown how, under appropriate hypotheses,
the property of (local) asymptotic regulation can be made robust
with respect to plant parameter variations, extending in this way
a celebrated property of linear regulators.

In the presence of plant parameter variations, the challenge is to
design a (parameter independent) controller in such a way that the
closed-loop system possesses a (possibly  parameter {\em
dependent})  attractive invariant set on which the regulated
output vanishes. The two issues of forcing the existence of such
an invariant set and of making the latter (locally or non-locally)
attractive are of course interlaced and this is precisely what, in
the past years, has determined the various scenarios under which
different solutions to the problem have been proposed. In the
paper \cite{HL}, for instance, this was achieved by assuming that
the set all feed-forward controls which force the regulated output
to be identically zero had to be generated by a single
(parameter-independent) {\em linear} system. This assumption was
weakened in \cite{CH3}, in \cite{DMI04} and, subsequently, in
\cite{BI03bis}, where it was replaced by the assumption that the
controls in question are generated by a single
(parameter-independent) {\em nonlinear} system, uniformly
observable in the sense of \cite{GK}.

The crucial observation that made the advances in \cite{BI03bis} and
\cite{DMI04} possible was the realization that the  two issues of forcing the
existence of an invariant set (on which the regulated variable vanishes) and of
making the latter attractive are intimately related to, and actually can be
cast as, the problem of designing a (nonlinear) {\em observer}. As a matter of
fact, the design method suggested in \cite{BI03bis} was based almost entirely
on the construction of a nonlinear ``high-gain" observer following the methods
of Gauthier-Kupka \cite{GK}, while the design method suggested in \cite{DMI04}
was based almost entirely on the construction of a nonlinear adaptive observer
following the methods of Bastin-Gevers \cite{BaGe} and Marino-Tomei
\cite{MaTo92}.

Having realized that the design of observers is instrumental in
the design of controllers which solve the problem in question, the
idea came to examine whether alternative options, in the design of
observers, could be of some help in weakening the assumptions even
further. This turns out to be the case, as shown in the present
paper, if the approach to the design of nonlinear observers
outlined by Kazantis-Kravaris (\cite{KaKr}) and then pursued by
Kresselmeier-Engel (\cite{Kress}), by Krener-Xiao (\cite{Krener})
and by Andrieu-Praly (\cite{AnPr}) is adopted.

While in all earlier contributions it was assumed that the
controls which force the regulated output to be identically zero
could be interpreted as outputs of a (in general nonlinear) system
having special observability properties (which eventually became
part of the controller), a crucial property highlighted in the
proof of Theorem 3 of \cite{AnPr} shows that no assumption of this
kind is actually needed.  The controls in question can {\em
always} be generated by means of a system of appropriate dimension
whose dynamics are linear but whose output map is a nonlinear
(and, in general, only continuous but not necessarily locally
Lipschitzian) map. Once this system is embedded in the controller,
boundedness of all closed-loop trajectories and convergence to the
desired invariant set can be guaranteed, as in the earlier
contributions \cite{BI03bis} and \cite{DMI04}, by a somewhat
standard paradigm which blends practical stabilization with a
small-gain property for feedback interconnection of systems which
are input-to-state stable (with restrictions).

The purpose of this paper is to provide a complete proof of how the results of
\cite{AnPr} can be exploited for the design of a controller solving the
problem  and also to show how some technical hypotheses used in the asymptotic
analysis of \cite{BIP} can be totally removed, yielding in this way a general
theory
 cast only on a very simple and
meaningful assumption. This paper is deliberately meant to present
only all the theoretical
 results needed to show the {\em existence} of the solution of
 the problem in question. Issues related to practical aspects involving
 constructive design and implementation will be dealt with in a forthcoming
 work.

The paper is organized as follows. In the next section the main
framework
 under which the problem is solved is presented and discussed. Then Section 3
 presents an outline of the main results concerning the existence of the output
 feedback regulator. Technical proofs of the results in this section are
 postponed to Appendices A and B. Section 4 concludes the paper
 with some
 with final remarks.\\[2mm]

\noindent {\bf Notation.}
 For $x \in \Real^n$, $|x|$ denotes the Euclidean norm and, for $\cal C$ a closed
 subset of $\Real^n$, $|x|_{\cal C}=\min_{y \in \cal C}|x-y|$ denotes the
 distance of $x$ from $\cal C$. For $\cal S$ a
subset of $\Real^n$, $\mbox{cl}\cal S$ and $\mbox{int}\cal S$ are
 the closure of $\cal S$ and the interior of $\cal S$
 respectively,
 and $\partial {\cal S}$ its boundary.
 For the smooth dynamical system $\dot x = f(x)$, the value at time $t$ of the solution passing
 through $x_0$ at time $t=0$ will be written as $x(t,x_0)$. Somewhere the more compact
 notation $x(t)$ will be used instead of $x(t,x_0)$, when the initial condition is
 clear from the context. A set $\cal S$ is said to be locally forward (backward)
 invariant for $\dot x = f(x)$ if there is a time $t_0>0$ ($-t_0<0$) such that, for
 each $x_0 \in \cal S$,  $x(t,x_0) \in \cal S$ for all $t \in [0,t_0)$
 ($t \in (-t_0,0]$). The set
 is locally invariant if it is locally backward and forward invariant. The set
 is (forward/backward) invariant if it is locally (forward/backward) invariant
 with $t_0=\infty$.
 For a locally Lipschitz function $V(t)$ we define the Dini's derivative of $V$ at $t$ as
 \[
  D^+ V(t) = \lim_{h \rightarrow 0^+} \sup {1 \over h} [V(t+h) - V(t)]\,.
 \]
 By extension, when $V(t)$ is obtained by evaluating $V$ along a
solution $x(t,x_0)$, we denote also
\beeq{\label{D^+}
D^+ V(x_0) = \lim_{h \rightarrow 0^+} \sup {1 \over h} [V(x(h,x_0)) -
V(x_0)]\,.
}
Note that if $\limsup=\lim$, this is simply $L_fV(x_0)$, the Lie derivative at
$x_0$ of $V$ along $f$.

 \section{The framework}

 \subsection{The problem of output stabilization and the main result} \label{sec21}
 We consider in what follows a nonlinear smooth system described by
 \beeq{\label{pla1}
  \ba{rcl} \dot z &=& f(z, \zeta)\\
  \dot \zeta &=& q(z,\zeta) + u\\[1mm]
 \ea
 }
 with state $(z,\zeta)\in \Real^n \times\Real$ and {\em control} input $u \in \Real$,
 and with unknown initial conditions $(z(0),\zeta(0))$
 ranging in a known arbitrary compact
 set $Z \times \Xi \subset \Real^n \times \Real$.

 Associated with (\ref{pla1}) there is a {\em controlled} output
 $e\in \Real$ expressed as \beeq{\label{error}
 e=h(z,\zeta)}
and a {\em measured} output $y \in \Real^{p}$  expressed as
\beeq{\label{output}
 y=k(z,\zeta)}
 in which $h: \Real^n \times \Real \to \Real$ and
  $k: \Real^n \times \Real \to \Real^p$  are smooth functions.

 For system (\ref{pla1})-(\ref{error})-(\ref{output}) the problem
 of semiglobal (with respect to $ Z\times \Xi$)
 { output stabilization} is defined as follows.
 Find, if possible,
 an output feedback controller of the form
 \beeq{\label{reg1}
 \ba{rcl}
 \dot \eta &=& \varphi(\eta, y) \\
 u &=& \varrho(\eta,y)
 \ea
 }
 with state $\eta\in \Real^\nu$ and a compact set $M \subset
 \Real^\nu$ such that, in the associated closed loop system
 \beeq{\label{closedloop}
 \ba{rcl}
 \dot z &=& f(z, \zeta)\\
  \dot \zeta &=& q(z,\zeta) + \varrho(\eta,k(z,\zeta))\\
  \dot \eta &=& \varphi(\eta, k(z,\zeta))\\
  e &=& h(z,\zeta)\,,\ea}
 the
 positive orbit of $Z\times \Xi \times M$ is bounded
  and, for each $(z(0),\zeta(0),\eta(0)) \in Z\times  \Xi  \times
  M$,
 \[\lim_{t \to \infty} e(t)=0\,.\]

 The problem at issue will be solved under the following main assumption
 which, roughly speaking, requires that  system
 \beeq{\label{zdadd}\ba{rcl}
 \dot z &=& f(z, \zeta)\ea}
viewed as a system with input $\zeta$ and output
 (\ref{output}), is ``stabilizable" by output feedback in an appropriate sense.
 In more precise terms the assumption in question is formulated as follows.\\[1mm]

  \noindent{\bf Assumption.} There exists a compact set
  ${\cal A} \subset \Real^{n}$, a smooth function $\alpha:
  \Real^n \to \Real$ and a smooth map
 $\Phi: \Real^p \to \Real$ such that:
    \begin{itemize}
     \item[(${\bf a}_1$)] the set ${\cal A}$ is locally asymptotically stable for the
     system
     \beeq{\label{zd}
        \ba{rcl}
          \dot z &=& f(z,\alpha(z))
        \ea
     }
     with a domain of attraction ${\cal D} \supset  Z$;
     \item[(${\bf a}_2$)] $h(z,\alpha(z))=0$ for all $z \in \cal A$.
     \item[(${\bf a}_3$)] $\Phi(k(z,\zeta)) = \zeta-\alpha(z)$
     for all $(z,\zeta) \in \Real^n \times \Real$.
    \end{itemize}
  Comments to this assumption are postponed after the next theorem
  which presents the main result of the paper.

 \begin{theorem} \label{MainTheorem}
 There exists an $m>0$, a controllable pair $(F,G) \in \Real^{m \times m}
 \times \Real^{m \times 1}$, a continuous function
 $\gamma: \Real^m \to \Real$ and, for any compact set $M \subset
 \Real^m$, a continuous function $\kappa: \Real^{p} \to \Real$,
  such that the controller
   \beeq{\label{reg}
   \ba{rcll}
   \dot \eta &=& F \eta + G u  &\qquad  \eta(0) \in M\\
   u &=& \gamma(\, \eta \,) + v & \\
   v &=& \kappa(y)
   \ea
  }
 solves the problem of semiglobal (with respect to $ Z\times \Xi$)
 output stabilization.
 \end{theorem}

 \begin{remark}
 It is worth noting that the previous Assumption could be  weakened by
 asking not for a static but rather for a dynamic output feedback stabilizer.
 More precisely  for all the results presented in the paper to hold, it suffices to
 assume the
 existence of a smooth system
 \[\ba{rcl}
 \dot \xi &=& \sigma(\xi, u_\xi)
 \\
 y_\xi&=&
\beta(\xi, u_\xi)\ea\] with initial state allowed to range on  a compact subset
$\Sigma$ of $\Real^l$, of a pair of smooth functions $\alpha: \Real^n \to
\Real$ and
 $\Phi: \Real^p \to \Real$, and of a compact set $\ca \subset \Real^n \times
 \Real^l$ such that:

 \medskip\noindent
 (${\bf a}_1^\prime$) the set $\ca$ is locally asymptotically stable for
 the system
 \[
     \ba{rcl}
          \dot z &=& f(z,\beta(\xi, \alpha(z)))\\
          \dot \xi &=& \sigma(\xi, \alpha(z))\,,
        \ea
 \]
 with a domain of attraction ${\mathcal D} \subset
 Z\times \Sigma$,

 \medskip\noindent(${\bf a}_2^\prime$)
 $h(z,\beta(\xi, \alpha(z)))=0$ for all $(z,\xi) \in \ca$,

 \medskip\noindent and condition (${\bf a}_3$) above holds.
This, however,
 is not deliberately pursued in what follows, as it  would only add unnecessary
 complications,
 without  any extra conceptual value. $\triangleleft$
   \end{remark}

 \subsection{Remarks on the framework and the result} \label{sec22}

Since a crucial requirement
 in
 the problem of semiglobal output stabilization is boundedness of the closed-loop
 trajectories,  it seems natural to postulate the existence of a  ``virtual"
 control $\zeta(\cdot)$ that keeps the trajectories of (\ref{zdadd}) bounded.
 Assumption ($\bf a_1$) does this by asking that the control in question be
 a state feedback law, namely a control
 $\zeta=\alpha(z)$ under which the trajectories of (\ref{zdadd}) asymptotically approach
 a compact set $\ca$. If this is the case, the
 change of
variables
\[
\chi = \zeta-\alpha(z)
\]
changes system (\ref{pla1}) into a system of the form
\beeq{\label{pla1bis}
  \ba{rcl} \dot z &=& f(z, \alpha(z)+\chi)\\
  \dot \chi &=& \tilde q(z,\chi) + u\\[1mm]
 \ea
 }
 in which
 \[\tilde q(z,\chi)=q(z,\alpha(z)+\chi)+\dst{\partial \alpha \over \partial z}f(z,
 \alpha(z)+\chi).
 \]
 Accordingly, the controlled output $e$ becomes
 \[
 e = h(z,\alpha(z)+\chi)\,.
 \]
 Clearly, to take advantage of the fact that the trajectories of
 (\ref{zd}) with initial conditions in $Z$, as required in part ($\bf a_1$) of the
 Assumption, are attracted by a
 compact set $\ca$, it would be desirable to have asymptotically $\chi$
 converging to zero.

 Of course since the same ``virtual" control
  appears in the map which expresses the controlled variable $e$,  and the latter is required to
  asymptotically decay to 0, it is also appropriate
 to assume,
 as done in  $({\bf a}_2)$, that $h(z,\alpha(z))$ vanishes
 on $\ca$. If this were to occur, in fact, then also
 the controlled variable $e$ would converge to zero and the problem would be
 solved.

 To make $\chi$ converging to zero, one might wish to
 appeal to (somewhat standard) ``high-gain" arguments and have
 $u=-k\chi$, which would be an admissible control law because, as required in part ($\bf
 a_3$) of the Assumption, $\chi=\Phi(y)$ is available from the measured
 output. However, it is well known (see e.g. \cite{TP}) that to
 have $\chi$ asymptotically converging to zero in a ``high-gain" scheme, it is
 somewhat necessary that the ``coupling" term $\tilde q(z,\chi)$
 between the upper and the lower subsystem of (\ref{pla1bis})
 asymptotically vanishes. More specifically, it is necessary that $\tilde
 q(z,0)$ be vanishing on the set $\ca$ to which the state $z$ of the
 upper subsystem converges if $\zeta$ decays to zero. Now, in
 general, there is no guarantee that $\tilde q(z,0)$ would vanish on $\ca$ and this is why a more elaborate
 controller has to be synthesized. As a matter of fact, the main
 result of the paper is that a suitable {\em dynamic} controller
 makes it sure that a property of this kind is achieved.

 A special case covered by the previous
 setup is the one in which system (\ref{zdadd}) can be given the form
 \beeq{\label{oldproblem}\ba{rcl}
 \dot z_1 &=& f_1(z_1)\\
\dot z_2 &=& f_2(z_1,z_2,\zeta)\,.\ea
} In this case, it is clear that the dynamics of $z_1$ is a
totally {\em autonomous} dynamics, which can be viewed as an
``exogenous" signal generator. This is the way in which the
classical problem of output regulation is usually cast. Depending
on the control scenario, the variable $z_1$ may assume different
 meanings. It may represent exogenous disturbances to be rejected
 and/or references to be tracked. It may also contain a set of
 (constant or time-varying) uncertain parameters affecting the controlled
 plant.

In this context, it is important to note that the proposed
framework encompasses  a number of problems which have been
recently
 addressed
 (see,
   among others, \cite{Kh2}, \cite{SIM},
  \cite{DMI04}, \cite{BI03}, \cite{BIM04}, \cite{CH3})  and rely
  upon various versions of the so-called
  ``minimum-phase" property. More specifically all the aforementioned works
   consider  systems having relative degree $r\ge 1$
 and normal form
  \beeq{\label{rdrpla}
  \ba{rcl}
  \dot w &=& s(w)\\
  \dot x_1 &=& f(w,x_1,C x_2)\\
  \dot x_2 &=& A x_2 + B \zeta\\
  \dot \zeta &=& q(w,x_1,x_2,\zeta) + u
  \ea
  }
  in which $w\in\Real^s$ represents an exogenous input,
  $x_1\in \Real^\ell$, $x_2\in \Real^{r-1}$, $\zeta\in \Real$ and $A,B,C$
  is a triplet in ``prime" form, with controlled and measured output
  respectively given by
  \[\ba{rcl}e &=& C x_2\\
  y &=& {\rm col}(x_2,\zeta)\,.\ea\]
  Note that ${\rm col}(x_2,\zeta) ={\rm col}(e, \dot e, \ldots,
e^{r-1})$. \footnote{The case in which $y = {\rm col}(x_2,\zeta)$
is known in the literature at the case of {\em partial state
feedback}, to emphasize the fact that not only the error $C x_{2}$
is available as measured output but also all its first $r-1$ time
derivatives. This is not a restriction, though, since -- as shown
for instance in \cite{EK} and \cite{TP} -- so long as convergence
from a compact set of initial conditions is sought, all components
of $y$ can always estimated by means of an ``approximate" observer
driven only by its first component $Cx_2$.}

  For this class of systems one of the main assumptions under which the problem of
  output regulation has been  solved is the one requiring that the ``zero
  dynamics"
  \beeq{ \label{zdsec2}
  \ba{rcl}
  \dot w &=& s(w)\\
  \dot x_1 &=& f(w,x_1,0)
  \ea
  }
  satisfy some stability requirement. For instance in \cite{SIM} the assumption in
  question asks for the  existence of a differentiable map $\pi: \Real^s \to
  \Real^{\ell}$  whose graph  $\ca' = \{(w,x_1) \in \Real^s \times
  \Real^{ \ell} \; : \; x_1 = \pi(w) \}$ is  invariant and locally
   exponentially stable for (\ref{zdsec2}), uniformly with respect to
   the exogenous variable $w$,  with a domain of attraction containing the
    assigned compact set
   of initial conditions. This assumption has been  substantially weakened in the recent
   work \cite{BI03} (see also \cite{BI03bis} and \cite{BIM04}) by asking a form of
   `` weak minimum-phase'' property, in which the compact set
   $\ca'$ above is rather the graph of a {\em set-valued} map.\footnote{ As shown in \cite{BI03} this assumption  is a straightforward
   consequence of the boundedness of the trajectories of (\ref{zdsec2}).}
   Under this  assumption, as shown in \cite{BIM04} (see also \cite{DMI04}), it is
   possible
   to argue the existence of a matrix $K$, designed via high-gain arguments,
   such that the system
 \beeq{\label{zd2sec2}
  \ba{rcl}
  \dot w &=& s(w)\\
  \dot x_1 &=& f(w,x_1,C x_{2})\\
  \dot x_2 &=& A x_2 + B K x_2\\
  \ea
  }
  possesses an {\em invariant} compact set $\ca=\{(w,x_1,x_2)\;:\; (w,x_1) \in \ca'\,,
   x_2 =0 \}$ which is locally exponentially stable for (\ref{zd2sec2}) with a domain
   of attraction containing the set of initial conditions.

  It is clear that the
    case
   described above perfectly fits in the  framework  presented in Section
   \ref{sec21},
   with the $(w,x_1,x_2)$ subsystem in (\ref{rdrpla}) playing the role of the
   $z$-subsystem in (\ref{pla1}), with the controlled  output and  measured
   output map in
   (\ref{error})-(\ref{output})
   respectively equal to $C x_{2}$ and  $(x_2,\zeta)$, and with assumptions
   ($\bf a_1$)-($\bf a_2$)-($\bf a_3$) automatically satisfied with
   $\ca=\{(w,x_1,x_2)\;:\; (w,x_1) \in \ca'\,, x_2 =0 \}$, $\alpha(z) =
   K x_2$ and $\Phi(y) = \qmx{-K & 1}y$.

   On the other hand, the framework set up in Section \ref{sec21}
   makes it possible to deal
   with some form of non-minimum-phase systems, extending in this way
   the class of systems which can be treated.

  \section{Main results}
  \subsection{ The basic approach}
  In this section we overview the main steps which will be followed to prove
  Theorem \ref{MainTheorem}. Technical proofs of the results given here
  are presented in Appendix \ref{SecProofProposition1}.

  We consider the closed-loop system (\ref{pla1}), (\ref{reg})
  which, after the change of coordinates
  \[
  \zeta \to \chi=\zeta -\alpha(z)\qquad \qquad \eta \to x=\eta - G \chi\,,
  \]
  can be rewritten as
 \beeq{\label{closedloopsys}
  \ba{rcl}
   \dot z &=& f_0(z) + f_1(z, \chi)\\[2mm]
   \dot x &=& F x  - G q_0(z)  - G q_1(z,\chi) + FG \chi\\[2mm]
   \dot \chi &=& q_0(z)  + q_1(z,\chi) + \gamma(x + G \chi) + v
   \ea}
 in which
 \beeq{\label{q0def}\ba{rcl}
 f_0(z) &:=& f(z,\alpha(z))\\[1mm]
 q_0(z) &:=& q(z,\alpha(z))  - \dst{\partial \alpha(z) \over \partial z} f(z,\alpha(z))
 \ea}
 and
 \[
 \ba{rcl}f_1(z,\chi) &:=&f(z,\alpha(z) + \chi) -f(z,\alpha(z) )\\[1mm]
\dst q_1(z,\chi) &:=& \dst q(z,\alpha(z)+ \chi) -q(z,\alpha(z) )  - {\partial
\alpha(z) \over \partial z} \left[ f(z,\alpha(z) + \chi) -f(z,\alpha(z)
)\right]\,.
  \ea\]
 Observe that we have $f_1(z,0) \equiv 0$ and $q_1(z,0) \equiv 0$ for all $z \in \Real^n$.

 In what follows, system (\ref{closedloopsys}) is seen
as a system with input $v$, output $\chi$ (which is, according to
 its definition and to assumption $({\bf a}_3)$, a nonlinear
 function of the original output $y$)
 and initial conditions contained in a
 set of the form
 $ Z \times X \times C$ in which
 $X\subset \Real^m$ and $C \subset \Real$ are compact sets dependent on $\Xi$
 and $M$.
 A controller of the form (\ref{reg}) solves the problem at issue
 if, for some map $\hat \kappa: \Real^n \to \Real^n$,
 the control law $v=\hat \kappa(\chi)$
 is such that  all
 trajectories of (\ref{closedloopsys}) originating from $ Z
 \times X \times C$  are bounded and
 \beeq{\label{chiwztozero}
 \lim_{t \to \infty} \chi(t) =0\,,
 \qquad
 \lim_{t \to \infty} |z(t))|_{\cal A}=0\,.
 }
 As a matter of fact, since
 systems (\ref{pla1})
 -(\ref{reg}) and
 (\ref{closedloopsys}) are diffeomorphic, boundedness of the trajectories of
 (\ref{closedloopsys}) with initial conditions in $
 Z \times X \times C$
 implies boundedness of the trajectories of
 (\ref{pla1})
 -(\ref{reg}) originating from
 $ Z \times \Xi \times M$. Furthermore,
 by virtue of assumption
 $({\bf a_2})$, since $h(\cdot)$ is a continuous function, condition (\ref{chiwztozero}) implies also that
 $\lim_{t \to \infty} e(t) =0$, namely that the
 problem of semiglobal output stabilization is solved.\footnote{Note
 that the map $\kappa(\cdot)$ to be used in (\ref{reg}) is
 actually $\hat \kappa(\Phi(\cdot))$.}

 By virtue of this fact, in the following we focus our attention on system (\ref{closedloopsys})
 and we prove that (\ref{closedloopsys}) controlled
 by $v=\hat \kappa(\chi)$
 has bounded trajectories and (\ref{chiwztozero}) hold.
 To this end, for notational convenience, denote
 \[
 p = \mbox{col}(z,x)
 \]
 and rewrite system (\ref{closedloopsys}) in the more compact form
  \beeq{\label{wpsys}
  \ba{rcl}
  \dot p &=& M(p) + N(p,\chi)\\[1mm]
  \dot \chi &=& H(p ) + K(p,\chi) + v
  \ea
 }
 in which
$M(\cdot)$ and $H(\cdot)$ are defined as
 \beeq{\label{Fdef}
  M(p)= \left ( \ba{c}
    f_0(z)\\[2mm]
    F x - G q_0(z)\\[2mm]
  \ea
  \right )
 }
 and
 \beeq{\label{Hdef}
  H(p)= q_0(z) + \gamma(x)
} and $N(\cdot)$ and
 $K(\cdot)$ are therefore suitable remainder functions which satisfy
   necessarily $N(p,0)=0$ and $K(p,0)=0$ for all $p$.
 Consistently set $P=Z \times X$ so that the initial conditions
 of (\ref{wpsys}) range in $ P \times C$. System (\ref{wpsys})
 is recognized to be a system in normal form
 with relative degree one (with respect to the input $v$ and output $\chi$)
 and zero dynamics given by
 \beeq{\label{zdwp}
 \ba{rcl}
  \dot p = M(p)\,.
 \ea
 }
 Thus, following consolidated knowledge about stabilization of {\em minimum-phase}
 nonlinear systems (see \cite{BI91}, \cite{Bacc}, \cite{TP}), the capability of
 stabilizing (\ref{wpsys})
 by output feedback is expected to strongly rely upon asymptotic properties
 of the zero dynamics (\ref{zdwp}). This is confirmed by the next two results
 showing that the existence of an asymptotically stable attractor for system (\ref{zdwp})
 is sufficient to achieve boundedness of trajectories and {\em practical} stabilization
 (Theorem \ref{theorempracticalreg}), which becomes {\em asymptotic}
 if the function $H(\cdot)$ vanishes on the attractor
 (Theorem \ref{theoremasymptoticreg}).
 These results, which in the context of this paper represent
 building blocks for proving Theorem \ref{MainTheorem}, are interesting by their
 own as they represent an  extension of well-known stabilization paradigms for
 systems with equilibria (see \cite{TP}) to the case of systems of the form
 (\ref{wpsys}),  with zero dynamics (\ref{zdwp}) possessing compact
 attractors.

 \begin{theorem} \label{theorempracticalreg}
 Consider system (\ref{wpsys}) with $M(\cdot)$ at least locally Lipschitz
 function and $N(\cdot)$, $H(\cdot)$, $K(\cdot)$ at least continuous
 functions. Let the initial conditions be in $ P \times
 C$. Assume that
 system (\ref{zdwp}) has a compact attractor $\cal B$ which is
 asymptotically stable with a domain
 of attraction ${\cal D} \supset P$. Then for all $\epsilon>0$
 there exists a $\kappa^\star>0$ such that for all $\kappa \geq \kappa^\star$
 the trajectories of (\ref{wpsys}) with $v= - \kappa \chi $
 are bounded and $\limsup_{t \to \infty} |\chi(t)| \leq \epsilon$ and
 $\limsup_{t \to \infty} |p(t)|_{\cal B} \leq \epsilon$.
 \end{theorem}

 \begin{theorem} \label{theoremasymptoticreg}
 In addition to the hypotheses of the previous theorem, assume
 that $\left . H(p) \right |_{\cal B} =0$. Then there exists a continuous function
 $\kappa: \Real \to \Real$ such that the trajectories of (\ref{wpsys}) with
 $v= \kappa(\chi)$ are bounded and $\lim_{t \to \infty} \chi(t) =0$ and
 $\lim_{t \to \infty} |p(t)|_{\cal B} =0$. If,
 additionally, $H(\cdot)$ and $K(\cdot)$ are locally Lipschitz and the set
 $\cal B$ is also locally exponentially stable for (\ref{zdwp}), then there
 exists $\kappa^\star>0$ such that for all $\kappa \geq \kappa^\star$ the
 same properties hold with $v= -\kappa \chi$.
 \end{theorem}

 For the proofs of these theorems the reader is referred to the sections
 \ref{Prooftheorempracticalreg} and \ref{Prooftheoremasymptoticreg} respectively.

 Motivated by these results (and in particular by Theorem
 \ref{theoremasymptoticreg}),  we turn our attention
 on the study of the zero dynamics (\ref{zdwp}) (with $M(\cdot)$ as in (\ref{Fdef}))
 and on the function $H(\cdot)$ in (\ref{Hdef}), by looking for the  existence of
 a  controller which guarantees the basic requirements behind Theorem
 \ref{theoremasymptoticreg}, with in particular $\left . H(p) \right |_{\cal B}=0$.
 Details in this direction are presented in the
 next subsection.

 \subsection{The properties of the ``core subsystem'' (\ref{zdwp})}
 The crucial result which will be proved in this part is that, under
 the assumption presented in section \ref{sec21}, there is a choice
 of the pair $(F,G)$ and of the map $\gamma(\cdot)$ which
 guarantee the existence of an asymptotically stable compact attractor $\cal B$
 for  (\ref{zdwp}), on which the function $H(\cdot)$ in (\ref{Hdef})
 vanishes. Moreover the projection of $\cal B$ on the $z$ coordinates coincides with
 $\cal A$. In view
 of the arguments discussed in the previous subsection,
this, along with an appropriate choice of $\kappa(\cdot)$ whose existence is
 claimed in Theorem \ref{theoremasymptoticreg},  substantially proves
 Theorem \ref{MainTheorem}.

 The result in question is proven in the next three
 propositions. To this end, note that  the core subsystem (\ref{zdwp})
 in the original coordinates $(z,x)$ is expressed as
 \beeq{\label{zxsys}
 \ba{rcl}
 \dot z &=& f_0(z)\\[1mm]
 \dot x &=& F x - G q_0(z)
 \ea
 }
 with initial condition in $Z \times X$.
 The first proposition is related to the
 first basic requirement behind Theorem \ref{theorempracticalreg}, namely
 the existence of a locally  asymptotically stable attractor for
 (\ref{zxsys}).

 More precisely, under the only requirement that $F$ is an Hurwitz matrix,
 it is shown the existence of  a set which is {\em forward invariant} and
  {\em locally asymptotically  stable} for (\ref{zxsys}).
 The set in question is described by the graph of a map.

 \begin{proposition} \label{propositiontau}
 Consider system (\ref{zxsys}) with  $\alpha(\cdot)$ such that
 condition (${\bf a}_1$) holds and let $(F,G)$ be any pair
 with $F$ Hurwitz. Then:

  \begin{itemize}
 \item[(i)]  there exists
 at least one continuous map $\tau: \Real^{n} \to \Real^m$
 such that the set
 \beeq{ \label{graphtaurest}
  \mbox{\rm graph}(\left . \tau \right |_{\ca}):=\{(z,x)\in \ca \times \Real^m \; :
  \; x=\tau(z)\}
 }
is forward invariant for (\ref{zxsys}).

  \item[(ii)] The set $\mbox{\rm graph}(\left . \tau \right |_{\ca})$
  is locally asymptotically stable for (\ref{zxsys})
  with a domain of attraction containing $Z \times X$. Furthermore the set
 in question is also locally exponentially stable for (\ref{zxsys}) if $\ca$ is
 such for (\ref{zd}).
 \end{itemize}
 \end{proposition}
 The proof of this proposition can be found in Section \ref{proofpropositiontau}.

{\begin{remark} Indeed, there might be many different continuous
maps $\tau$ having the property (i) of proposition 1. However, it
turns out that if $\ca_0$ is any compact {\em subset} of  $\ca$
which is {\em invariant} for (\ref{zd}), then for each $z\in\ca_0$
there is one and only one $x_z\in \Real^m$ such that the set
$\cup_{z\in\ca_0}\{(z,x_z)\}$ is invariant for (\ref{zxsys}). In
particular,
\[
x_z=-\int_{-\infty}^0
 e^{-Fs}G{q_0}(z(s,z)) ds
\]
where $z(s,z)$ denotes the value at time $t=s$ of the solution of
$\dot z=f_0(z)$ passing through $z\in
\ca_0$ at time $t=0$ (see \cite{BI04}). $\triangleleft$
\end{remark}}

 The second crucial requirement imposed by Theorem \ref{theoremasymptoticreg} is
that the function $H(\cdot)$ in (\ref{Hdef}) vanishes on the
asymptotically stable attractor $\mbox{\rm graph}(\left . \tau
\right |_{\ca})$. Here is where the precise choice of the pair
$(F,G)$ and of the map $\gamma(\cdot)$ play a role. In particular
note that, by definition of $H(\cdot)$ in (\ref{Hdef}) and of
$\mbox{\rm graph}(\left . \tau \right |_{\ca})$ in
(\ref{graphtaurest}), it turns out that
 \beeq{\label{HBpAp}
 \left . H(p) \right |_{\mbox{\rm graph}(\left . \tau \right |_{\ca})}=
 \left . \left ( {q_0}(z) +
 \gamma \circ \tau (z)\right )\right |_{\ca}
 }
 from which it is apparent that $\gamma(\cdot)$ should be chosen to satisfy
 $\gamma \circ \tau (z)=-{q_0}(z) $ for all $z$ in $\ca$.
 It is easy to realize that the possibility of choosing $\gamma(\cdot)$ in this way
 is intimately related to the fact that the map $\tau$ satisfies the {\em partial
(with respect to ${q_0}(\cdot)$) injectivity condition}
 \beeq{\label{inj1}
  \tau(z_1) = \tau(z_2)
  \qquad \Rightarrow \qquad
  {q_0}(z_1) = {q_0}(z_2) \qquad \mbox{for all } z_1, z_2 \in {\cal
  A}\,.
 }
 As $\tau$  is dependent on the pair $(F,G)$, the next natural point to be addressed is
 if there exists a choice of $(F,G)$ yielding the desired property for $\tau(\cdot)$.
 To this end is devoted the next proposition which claims that, indeed, there exists a
 suitable choice of $(F,G)$, with $F$ Hurwitz, such that the associated map $\tau(\cdot)$
 satisfies the required partial injectivity condition.
 Besides others technical constraints on
 the choice of $F$ which will be better detailed in the proof of the Proposition
 \ref{propositionvarrho}, the main requirement on $F$ is given by
 its dimension which is required to be sufficiently large with respect to the
 dimension of $z$.
 \begin{proposition} \label{propositionvarrho}
 Set
 \[
 m = 2 + 2 n\,.
 \]
 Then there exist a controllable pair $(F,G) \in
 \Real^{m \times m} \times \Real^{m\times 1}$, with $F$ a Hurwitz matrix,
 and a  class-$\cal K$ function $\varrho:\Real^+ \to \Real^+$
 such that
 \beeq{\label{partialinjection}
 |{q_0}(z_1) - {q_0}(z_2) | \leq \varrho(|\tau(z_1) - \tau(z_2)|)
 \qquad \mbox{for all } z_1, z_2 \in \ca
 }
 in which $\tau(\cdot)$ is a map (associated to $F$) with the properties indicated
 in Proposition \ref{propositiontau}.
 \end{proposition}
 For the proof of this proposition the reader is referred to Section
 \ref{proofpropositionvarrho}.

 \begin{remark} \label{remarkFG}
  Going through the proof of the previous proposition, it turns out that the pair
  $(F,G)$ be chosen as any $(2 n +2)$-dimensional real representation
  of the $(n+1)$-dimensional complex pair $(F_c,G_c)$, with $F_c=\mbox{diag}(\lambda_1, \ldots, \lambda_{r+1})$, $G_c=(g_1, \ldots,
  g_{r+1})\tr$ in which $g_i$
  are arbitrary not zero real numbers and  $\lambda_i$ are $n+1$ complex
  numbers taken arbitrarily outside a set of zero Lebesgue measure and with
  real part smaller than $\ell$, a real number related to the Lipschitz constant
  of $f_0(\cdot)$ (see Proposition \ref{Proptauc1}).
  $\triangleleft$
 \end{remark}

  It turns out that the injectivity property (\ref{partialinjection}) is
  a sufficient condition for the map $\gamma(\cdot)$ to exist.
  This is formalized in the next final proposition, proved in Section
  \ref{proofpropositionfinal}, which states that if
 (\ref{partialinjection}) holds then there exists a map $\gamma(\cdot)$ which makes
 $H(\cdot)$ vanishing on the attractor $\mbox{\rm graph}(\left . \tau \right |_{\ca})$.
 The map $\gamma(\cdot)$ can be claimed, in general, to be only continuous. It is also
 Lipschitz in the special  case in which the class-${\cal K}$ function $\varrho(\cdot)$ in
 (\ref{partialinjection}) is such.

 \begin{proposition}\label{Propositionfinal}
 Let $\tau(\cdot)$ be a continuous map satisfying (\ref{partialinjection}) with
 $\ca$ a closed set.
 Then there exist a continuous map $\gamma: \Real^m \to \Real$
 such that
 \beeq{ \label{qgammatau}
 {q_0}(z) +
 \gamma \circ \tau (z) =0
\qquad \forall \; z \in {\ca}\,.
 }
 If, in addition, the function $\varrho(\cdot)$ in (\ref{partialinjection}) is
 linearly bounded at the origin, then the map $\gamma$ is Lipschitz.
 \end{proposition}

  Combining the results of all the previous propositions, it appears that
  it is sufficient to choose the pair $(F,G)$ of suitable
  dimension (with $F$ Hurwitz) according to Proposition \ref{propositionvarrho} and
  to   choose $\gamma(\cdot)$ in order to satisfy relation (\ref{qgammatau}).
  In fact, so doing, we are guaranteed that the compact set
  ${\cal B}=\mbox{graph}(\left . \tau \right |_\ca)$ is locally asymptotically
  stable for (\ref{zxsys}) with the map (\ref{Hdef}) which is vanishing on $\cal B$.
  This, indeed, makes it possible to apply Theorem \ref{theoremasymptoticreg} and
  to conclude the existence of a continuous function $\kappa(\cdot)$,
    completing in this way  the synthesis of the controller.

 \begin{remark}
 The reader who is  familiar with  recent developments of the theory
 of
  nonlinear state observers will
 find interesting to  compare the previous results
  with the design method
 proposed by Kazantzis  and Kravaris  in \cite{KaKr}
 and pursued in  \cite{Kress}, \cite{Krener} and \cite{AnPr}.
 In the
 framework of \cite{KaKr}, system (\ref{zxsys}) can be identified with the
 cascade of an ``observed"
 system $\dot z = f_0(z)$ with output $y_z =q_0(z)$
 driving an ``observer" $\dot x = F x - G y_z$. If the map
 $\tau(\cdot)$ has a left inverse $\tau^{-1}_{\ell}(\cdot)$, the
 observer in question provides a state
  estimate $\hat z = \tau^{-1}_{\ell}(x)$.
 Such a left-inverse, as shown in \cite{AnPr},
 always exists provided that the dimension of $x$ is sufficiently
 large,
  if the pair $(f_0, q_0)$ has
 appropriate {\em observability} properties. In the present context of output
 stabilization, though, left invertibility of $\tau(\cdot)$
 is not needed.
 In fact, what the  controller is expected to do is
only the reproduction of the output $q_0(z(t))$ and not of the full state
$z(t)$ of the ``observed system". This motivates the absence of observability
 hypotheses on
 the pair $(f_0, q_0)$.
 $\triangleleft$
 \end{remark}

\section{Conclusions}
 This paper is focused on the existence of an output feedback
law that asymptotically steers to zero a given controlled variable, while
keeping all state variables bounded, for any initial conditions in a fixed
compact set. The proposed framework encompasses and extends a number of
existing results in the fields of output feedback stabilization and output
regulation of nonlinear systems. The main assumption under which the  theory is
developed is the existence of a state feedback control law able to achieve
boundedness of the trajectories of the zero dynamics of the controlled plant.
In this sense the result presented here is applicable for a wide class of
non-minimum-phase nonlinear systems not tractable in existing frameworks. In
the paper only results regarding the existence of the controller solving the
problem at hand have been presented while practical aspects involving its
design and implementation are left to a forthcoming work.

\appendix

\section{Converse Lyapunov result} \label{AppendixA}
Consider a system of the form
 \beeq{
 \label{uno}
 \dot p = f(p) \qquad p \in \Real^n \,,
 }
 in which $f(p)$  is a  $C^k$ (with $k$ sufficiently large)
 function, with initial condition ranging over a fixed {\em compact} set $P$.
 For system (\ref{uno}) assume the existence of a compact
 set ${\cal B} \subset \Real^n$ which
 is forward invariant and asymptotically stable for
 (\ref{uno}),
 with a domain of attraction ${\cal D} \supset P$.
 More precisely, by setting
 \[
 |p |_{{\cal B} / \cal D} = \left ( 1 + {1 \over
 |p|_{\partial \, \mbox{\small \rm cl}\, \cal D}}\right ) |p|_{\cal B}\,,
 \]
 we assume that the set $\cal B$ satisfies the following two properties:

 \noindent{\em Uniform stability:} there exists a class ${\cal K}$ function
$\varphi $ such that for any $\alpha>0$
 \[
 |p_0|_{\cb / \cal D} \leq \alpha \qquad \Rightarrow \qquad
 |p(t,p_0)|_{\cb / \cal D} \leq \varphi(\alpha)
 \quad \forall \;t \geq 0\,;
 \]
 \noindent{\em Uniform {attractivity}:} there exists a continuous function
 $T: \Real^+ \times \Real^+ \to \Real^+$ such that for any $\alpha>0$ and
 $\epsilon>0$
 \[
 |p_0|_{\cb / \cal D}
  \leq \alpha \qquad \Rightarrow \qquad
 |p(t,p_0)|_{\cb / \cal D} \leq \epsilon \quad
 \forall \; t \geq T(\alpha, \epsilon)\,.
 \]

 We say that $\cb$ is also {\em locally exponentially stable} for (\ref{uno})
 if there exist $M\geq 1$, $\lambda>0$ and $c_0>0$ such that
 \[
 |p_0|_{\cb / \cal D}  \leq c_0 \qquad \Rightarrow
 \qquad |p(t,p_0)|_{\cb / \cal D} \leq M e^{-\lambda t} |p_0|_{\cb / \cal D}\,.
 \]

In this framework it is possible to formulate the following
converse Lyapunov result which claims the existence of a {\em
locally Lipschitz} Lyapunov function vanishing on the attractor.
The result is not formally proved as it can be easily deduced by
the arguments presented in \cite{Yoshizawa} (see in particular
Theorem 22.5 and the related Theorems 22.1 and 19.2 in the quoted
reference).

 \begin{theorem}\label{TheoremYoshisawa}
 Under the above uniform stability and uniform attractivity
conditions, there exists a {continuous} function $V: {\cal D} \to \Real$
 with the following properties:
 \begin{itemize}
 \item[(a)] there exist class ${\cal K}_\infty$ functions $\underline a(\cdot)$,
 $\overline a(\cdot)$ such that
 \[
 \underline a(|p|_{\cb / \cal D}) \leq V(p) \leq \overline a(|p|_{\cb / \cal
 D})\qquad {\forall\;  p \in \cal D}\,;
 \]
 \item[(b)] there exists $c>0$ such that
\[
 D^+V(p) \leq -c V(p) \qquad \forall\; p \in {\cal D}\,;
\]
 \item[(c)] for all $\alpha>0$ there exists $L_\alpha>0$ such that
 for all $p_1$, $p_2 \in \cal D$ such that $|p_1|_{\cb / \cal D} \leq \alpha$,
 $|p_2|_{\cb / \cal D} \leq \alpha$ the following holds
 \[
    |V(p_1) - V(p_2)| \leq L_\alpha |p_1 - p_2|\,.
 \]\
 \end{itemize}
 If $\cal B$ is also locally exponentially stable for (\ref{uno}) then
 property (a) holds with $\underline a(\cdot)$,  $\overline a(\cdot)$ linear
 near the origin.
   \end{theorem}

With this result at hand, it is also possible to formulate a local
Input-to-State Stability result for system (\ref{uno}) forced by
an external signal. This is formalized in the next lemma.

 \begin{lemma} \label{lemma2AppB}
 Let $x:\Real_+\rightarrow \Real^m$ be a $C^0$ function.
Consider the system
 \beeq{\label{32zd}
 \dot p = f(p) + \ell(p,x(t))
 }
 in which $p \in \Real^n$ and $\ell(p,0)=0$  for all $p\in \Real^n$.
%
 The functions $f(\cdot), \ell(\cdot)$ are $C^1$. Suppose that system (\ref{uno})
 satisfies the assumptions expressed before.
 Then there exist functions $\beta(\cdot,\cdot)$ and $\gamma(\cdot)$,
 respectively of class $\cal KL$ and $\cal K$, and a $d^\ast>0$ such that if
 \beeq{\label{apprest}
 |p_0|_{\cal B} \leq d^{\,\ast} \qquad \mbox{and}
 \qquad |x(t)| \leq d^{\,\ast} \quad \mbox{for all }  t \geq 0
 }
  then the right maximal interval of definition of $p(t,p_0)$ is $[0,+\infty
)$ and we have
 \beeq{ \label{boundappb}
|p(t,p_0)|_{\cal B}
 \leq \max\left\{\beta(|p_0|_{\cal B},t)
\,  ,\, \gamma (\dst \max_{\tau \in [0,t]}
 |x(\tau)|)
\right\}
  \qquad \mbox{for all }  t \geq 0.
  }
  If the set $\cal B$ is also locally exponentially stable for (\ref{uno})
  then there exist  $N>1$, $k >0$ and $\bar \gamma>0$ such that (\ref{boundappb}) modifies as
 \beeq{ \label{expboundappb}
 |p(t,p_0)|_{\cal B}
 \leq N e^{- k t} |p_0|_{\cal B}+ \bar \gamma \dst \max_{\tau \in [0,t]}
 |x(\tau)|
  \qquad \mbox{for all }  t \geq 0.
  }
 \end{lemma}

 \begin{proof}
 Pick $\beta>0$ such that if $|p|_\cb \leq \beta$ then $p \in \cal D$
 and note that there exists $d_\beta>1$ such that for all $p$
 satisfying $|p|_\cb  \leq \beta$ then
 \beeq{\label{compdist}
  |p|_{\cb / \cal D} \leq d_{\beta} |p|_\cb.
 }
 As $\ell(\cdot)$ is
 differentiable and $\ell(p,0)=0$, there is an $\hat \ell>0$ such that
 for all $|p|_\cb \leq \beta$ and $|x| \leq d^\ast$
 \[
 |\ell(p,{x})| \leq \hat \ell |{x}| \,.
 \]

So consider the Lyapunov function $V$ given by Theorem \ref{TheoremYoshisawa}.
By using properties  (b) and (c) of this theorem setting $\alpha = d_\beta
\beta$, we obtain for the system (\ref{32zd}), so long as $|p_1|_\cb < \beta$
and $|x| \leq d^\ast$,
\beeq{
\label{Dpsteps}
\begin{array}[b]{rcl} D^{\,+}V(p_1,x) &=& \dst
\limsup_{h\to 0^+}{1\over h}[V(p(h,p_1))-V(p_1)]\\[2mm]
&= &
\dst \limsup_{h\to 0^+}{1\over h} [V(p_1 + h {{f}}(p_1) + h \ell
(p,x)) - V(p_1)]\\[2mm]
&\leq & \dst \limsup_{h\to 0^+}{1\over h} [V(p_1 + h {{f}}(p_1) + h
\ell (p_1,x))
- V(p_1 + h {{f}}(p_1))] \\[2mm]
&& \qquad \qquad + \dst \limsup_{h\to 0^+}{1\over h} [ V(p_1 + h {{f}}(p_1)) - V(p_1) ]\\[1mm]
&\leq& \dst \limsup_{h\to 0^+}{1\over h} L_{\alpha}  h \ell(p_1,x) -c
V(p_1) \leq L_{\alpha} \hat \ell |x| - c V(p_1) \,.
\end{array}
}
Now  assume (\ref{apprest}) holds. Let $p(t,p_0)$, shortly rewritten
$p(t)$, be the corresponding solution of (\ref{32zd}). Let $[0,T_0)$
be its right maximal interval of definition when restricted to take
values in the open set $\{p: |p|_\cb < \beta\}$.
(\ref{Dpsteps}) holds for $p_1=p(t)$ and all $t$ in $[0,T_0)$.
This implies
\beeq{\label{LP1} V(p(t)) \le
 e^{-c(t-t_0)}V(p_0) +{L_\alpha\hat \ell \over c}\max_{\tau\in
 [0,t]}|x(\tau)|
\qquad \forall t\in [0,T_0)\,. } This, in view of property $(a)$ in Theorem
 \ref{TheoremYoshisawa}, yields
 \beeq{\label{gain1}
 |p(t)|_{\cb } \le |p(t)|_{\cb /\cal D} \leq
 \underline a^{-1} (2 e^{-ct} \, \overline a(|p_0|_{\cb / \cal D}))
 +
 \underline a^{-1} (\dst 2{L_\alpha \hat \ell \over c} \max_{\tau\in
 [0,t]}|x(\tau)|)
\qquad \forall t\in [0,T_0)\,.
 }
By using (\ref{compdist}), it follows that if
 $d^{\,\ast}$ is chosen so that
 \[
 d^\star \leq \min \{{c \over 2 L_\alpha \hat \ell} \, \underline
    a({\beta / 3})\,,
 {1 \over d_\beta} \, \overline a^{-1}({1 \over 2} \, \underline
    a({\beta / 3}))
 \}
 \]
we have
\[ |p(t)|_{\cb } < \beta \qquad \forall t\in [0,T_0) \,  .\]
From the definition of $T_0$, it must be infinite. So we have established that
(\ref{gain1})
 holds for all $t\geq 0$ if (\ref{apprest}) is satisfied. This
 proves the first part of the result. The second part of the result,
 namely that under exponential stability the bound (\ref{expboundappb}) holds,
 follows immediately by (\ref{gain1}) by using the fact that the functions $\underline
 a(\cdot)$ and $\overline a(\cdot)$ can be linear near the origin.
 $\triangleleft$
\end{proof}


 \section{Proofs}\label{SecProofProposition1}

 \subsection{Proof of Theorem \ref{theorempracticalreg}}
 \label{Prooftheorempracticalreg}

 \begin{proof}
 The study of the feedback interconnection (\ref{wpsys}) can be done
 by means of arguments which are quite similar to those
 used in \cite{IsBook2} to prove some of the main stabilization results
 of \cite{TP}. In doing this, we take advantage of the
 converse Theorem \ref{TheoremYoshisawa} presented in Appendix \ref{AppendixA}.

Let $V:{\cal D} \to \Real$ be the function given by
Theorem \ref{TheoremYoshisawa}.
Pick a number $a>0$ such that $C \subset
B_a :=\{\chi  \in \Real  \; : \; |\chi | \leq a\}$ and $P \subset V^{-1}([0,a])$
 (which is possible because of property $(a)$ in Theorem
\ref{TheoremYoshisawa}). Define
\[
\hat c = \max_{(p,\chi)\in V^{-1}([0,a+1])\times B_{a+1}}|H(p) + K(p,\chi)|\,.
\]
Also, since $N(p,\chi)$ is differentiable and vanishes at
$\chi=0$, there is a number $\hat n$ such that
 \[
 |N(p,\chi)|\le \hat n|\chi|
 \qquad \forall (p,\chi) \in V^{-1}([0,a+1]) \times B_{a+1}.
 \]
Finally, by property $(c)$ in
Theorem \ref{TheoremYoshisawa}, there is a number $L_V$ such that
\[
|V(p_1)-V(p_2)| \le L_V|p_1-p_2|\qquad \forall (p_1,p_2) \in
V^{-1}([0,a+1])^2.
\]

Then, by choosing  $v = - k \chi$ in the $\chi$-dynamics in
(\ref{wpsys}), we get (see notation (\ref{D^+}))
\beeq{\label{dpchi}
D^+ |\chi| \leq -\kappa |\chi|  + \hat c
\qquad \forall (p,\chi)\in V^{-1}([0,a+1])\times B_{a+1}
}
Also, by following the same lines as in (\ref{Dpsteps}), we get
\beeq{\label{Vtrick1}  D^{\,+}V(p)
 \leq  L_V \hat n|\chi|- c V(p)
 \qquad \forall (p,\chi) \in V^{-1}([0,a+1]) \times B_{a+1}.
.
}

So now consider a solution $(p(t),\chi (t))$ issued from a point in
$P\times C \subset V^{-1}([0,a]) \times B_a$. Let $[0,T_1)$ be its
right maximal interval of definition when restricted to take values
in the open set int$\left(V^{-1}([0,a+1]) \times B_{a+1}\right)$. It follows that
both (\ref{dpchi}) and (\ref{Vtrick1}) hold for $(p(t),\chi (t))$
when $t$ is in $[0,T_1)$. They give successively, for all $t$ in $[0,T_1)$
\begin{eqnarray*}
|\chi (t)|&\leq &e^{-\kappa t}\,  a + \frac{\hat c}{\kappa }\left(1-e^{-\kappa t}\right)
\\
V(p(t))&\leq &
e^{-ct}\,  a  + L_V \hat n \left(\frac{\hat c}{\kappa }\frac{1-e^{-ct}}{c}
+\frac{e^{-c t}-e^{-\kappa  t}}{\kappa -c}\left[a-\frac{\hat c}{\kappa }\right]\right)
\\&\leq &
a  + L_V \hat n\left(\frac{\hat
c}{c\kappa }+\frac{a}{\kappa -c}\right)
\end{eqnarray*}
Hence, by selecting $\kappa $ to satisfy
\[
\kappa >\max\left\{2\hat c , (c+3aL_V \hat n) ,
\frac{3L_V \hat n \hat c}{c}\right\},
\]
we get, for all $t$ in $[0,T_1)$
\\[0.7em]\null \hfill $\displaystyle
|\chi (t)| \leq  a + \frac{1}{2}
$\hfill and \hfill $\displaystyle
V(p(t)) \leq  a + \frac{2}{3}
$\hfill \null \\[0.7em]
This says that the solution
remains in $V^{-1}([0,a+\frac{2}{3}]) \times B_{a+\frac{1}{2}}$. So
from its definition, $T_1$ is infinite.
Then, from (\ref{dpchi}), we get
\[ \limsup_{t\rightarrow +\infty } |\chi (t)| \leq \frac{\hat c}{\kappa } \]
With (\ref{Vtrick1}), this in turn implies
\[ \limsup_{t\rightarrow + \infty } V(p(t)) \leq \frac{L_V \hat n
\hat c}{\kappa } .
\]
In view of property $(a)$ in Theorem
 \ref{TheoremYoshisawa}, the latter yields
\[ \limsup_{t\rightarrow + \infty } |p(t)|_{\cb } \leq
 \underline a^{-1}\left( \frac{L_V \hat n \hat c}{\kappa } \right ).
\]
So the conclusion of Theorem \ref{theorempracticalreg} holds if we
further impose to $\kappa $ to satisfy
\\[0.7em]\null \hfill $\displaystyle
\kappa >\max\left\{
\frac{\hat c}{\epsilon} ,
\frac{L_V \hat n \hat c}{ \underline a(\epsilon)}
\right\}
$\hfill
$\triangleleft$

\end{proof}

 \subsection{Proof of Theorem \ref{theoremasymptoticreg}}  \label{Prooftheoremasymptoticreg}

 \begin{proof}
 The proof of this result follows by standard small gain arguments.
  Let $\kappa(\chi) = -\alpha(\chi)$ where $\alpha(\cdot)$ is a
 continuous function such that $\alpha(0)=0$ and $\chi\alpha(\chi)>0$
 $\forall \, \chi \neq 0$.
 By mimicking the proof of Theorem \ref{theorempracticalreg}
 it is possible to show that for any $\epsilon>0$ there exist a $\kappa^\star>0$
and a $T>0$
such that, if $|\alpha(|\chi|)|
 \geq \kappa^\star |\chi|$ then each trajectory of the closed-loop system
issued from the compact set $P\times C$ satisfies
\\[0.7em]\null \hfill
$\displaystyle |p(t)|_\cb \leq 2\epsilon$\hfill and \hfill
$\displaystyle |\chi(t)| \leq 2\epsilon$\hfill  for all $t\geq T$\hfill \null .
\par\vspace{0.5em}

Observe that Lemma \ref{lemma2AppB} applies to the $p$-component
of the closed loop solution. So let $d^\star$ be given by this lemma. By picking
$\epsilon$ above satisfying $2\epsilon \leq d^\star$ and by applying Lemma \ref{lemma2AppB},
(after time $T$,) we obtain
 \beeq{ \label{boundappbnew}
|p(t)|_{\cal B} \leq \max\left\{\beta(|p(T)|_{\cal B},t-T) \,  ,\, \gamma (\dst
\max_{\tau \in [T,t]} |\chi(\tau)|) \right\}
  \qquad \mbox{for all }  t \geq T \,.
  }

With the properties of the functions $H$ and $K$, there
exist of class-$\cal K$ functions $\varrho_h(\cdot)$ and $\varrho_q(\cdot)$
such that
  \[
  |H(p)| \leq \varrho_h(|p|_{\cal B}) \qquad
  |K(p,\chi)| \leq \varrho_k(|\chi|)\,.
 \]
 Clearly $\varrho_h(\cdot)$ and $\varrho_k(\cdot)$ can be taken
linearly bounded at the origin if $H(\cdot)$ and $K(\cdot)$ are locally Lipschitz.
We obtain, for all $(p,\chi )$,
$$
D^+|\chi |\; \leq \; \varrho_h(|p|_{\cal B})  + \varrho_k(|\chi|) -|\alpha (\chi)|\,.
$$
So let us choose $\alpha(\cdot)$ so that
 \[
 |\alpha(\chi) | \geq 3 \max \{\varrho_h(\bar \gamma^{-1}(|\chi|)), \,
 \varrho_k(|\chi|), \,
 \kappa^\star |\chi| \} + |\chi |
 \]
where $\bar \gamma(\cdot)$ is a class $\cal K$ function such that
 $\bar \gamma \circ \gamma (s)<s$ for all $s \in \Real^+$
with $\gamma $ given by Lemma 1 (see (\ref{boundappb})). This
gives
$$
D^+| \chi |\; \leq \; -|\chi |+ \left[
\varrho_h(|p|_{\cal B}) - \varrho_h(\bar \gamma^{-1}(|\chi|)) \right]\,  .
$$
So for the closed loop solution, we get
$$
|\chi (t)|\; \leq \; \max\left\{ \exp(-(t-T)) |\chi (T)| \,  ,\, \sup_{s\in [T,
t)}\bar \gamma (|p(s)|_{\cal B}) \right\}
$$
for all $t \geq T$.
 From this and (\ref{boundappbnew}) the first claim of the Theorem follows by small
 gain arguments. The second claim of the theorem
immediately follows by the previous arguments and by
(\ref{expboundappb}) in Lemma \ref{lemma2AppB}.$\triangleleft$
 \end{proof}

 \subsection{Proof of Proposition \ref{propositiontau}} \label{proofpropositiontau}

 Let ${\cal O}(Z)$ {denote} the positive orbit of $Z$
 under the flow of
\beeq{\label{zdbold}
 \dot z = f_0(z)\,,
 }
 which is a bounded and forward invariant set for (\ref{zdbold}), such that $\ca \subset
 {\cal O}({Z})$. Moreover let $\hat {\cal O}({Z})$ be a compact strict
 superset of ${\cal O}({Z})$ such that $\hat {\cal O}({Z}) \subset \cal D$ and
 define the system
 \beeq{ \label{hatbzsys}
 \dot {\hat z} = a_0(\hat z){{f}_0}(\hat z)
 }
 in which $a(\hat z): \Real^{n} \to \Real$ is any bounded smooth function
 such that
  \[
  a_0(\hat z) = \left \{\ba{ll}
  1 & \qquad \hat z \in {\cal O}({Z})\\
  0 & \qquad \hat z \in \Real^{n} \setminus \hat {\cal O}({Z})\,.
  \ea
   \right .
  \]
  Let $\hat z(t,z_0)$ and $z(t,z_0)$ denote the flows
  of (\ref{hatbzsys}) and, respectively, (\ref{zdbold}) and note that, as
  a consequence of the fact that ${\cal O}({Z})$ is forward
  invariant and that systems (\ref{zdbold}) and (\ref{hatbzsys}) agree on
  ${\cal O}({Z})$, it turns out that
  \beeq{\label{zzhatagree}
  \hat z(t,z_0) = z(t,z_0) \qquad \mbox{for all } z_0 \in {\cal O}({
  Z}) \; \mbox{and }
  t \geq 0\,.
  }
 Moreover note that for any $\hat z_0 \in \Real^{n}$, (\ref{hatbzsys}) has a
 unique  solution $\hat z(t,\hat z_0)$ which is defined and bounded on $ t
 \in (-\infty,  \infty)$.

  Define now
  \beeq{\label{taumapdef}
  \ba{rcccl}
 \tau &:& \Real^{n} &\to& \Real^{m}\\
 && {z} & \mapsto & \dst \int_{-\infty}^0
 e^{-Fs}G{{q}_0}(\hat {z}(s,{z})) ds
 \ea
  }
  which, as a consequence of the fact that $F$ is Hurwitz and
  ${{q}_0}(\hat {z}(s,{z}))$ is bounded and continuous in $z$
  for any $s\in \Real$, is a well-defined {\em continuous} map.  We show now
  that graph$(\left . \tau \right |_\ca)=\{(z , \xi) \in \ca \times \Real^m \; : \;
  x= \tau(z)\}$ is a forward invariant set for (\ref{zxsys}).
  For, pick $z_0 \in \ca$ and $x_0\in \Real^m$, let $(z(t,z_0), x(t,z_0,x_0))$
  denote the value at time $t$ of the solution of (\ref{zxsys}) passing through
  $(z_0,x_0)$ at time $t=0$, and note that for all $t \geq 0$ (using (\ref{zzhatagree}))
  \beeq{\label{invproof}
  \ba{rcl}
  x(t,z_0,\tau(z_0)) &=& \dst e^{Ft} \tau(z_0) +  \int_{0}^t e^{F(t-s)} G
  {{q}_0}(z(s,z_0)) ds\\[2mm]
  &=& \dst e^{Ft} \int_{-\infty}^0 e^{-Fs} G {{q}_0}(\hat z(s,z_0))ds
  +  \int_{0}^t e^{F(t-s)} G
  {{q}_0}(z(s,z_0)) ds\\[2mm]
  &=& \dst \int_{-\infty}^t e^{F(t-s)} G
  {{q}_0}(\hat z(s,z_0)) ds = \int_{-\infty}^0 e^{-F s} G {{q}_0}(
   \hat z(s+t,z_0)) ds\\[2mm] &=& \tau(\hat z(t,z_0)) =
   \tau(z(t,z_0))\,.
  \ea
  }
  This, along with the fact that $\ca$ is forward invariant for (\ref{zdbold}) and
  subset of ${\cal O}(Z)$,
  proves that graph$(\left . \tau \right |_\ca)$ is {\em forward invariant}
  for (\ref{zxsys}).

 We prove now item (ii) of the proposition. To this end note that, by
 (\ref{invproof}),
 it follows that
 \[
 L_{a_0f_0}\tau(z) = F \tau(z) - G {{q}_0}(z)
 \]
 for all $z \in {\ca}$. Defining $\tilde x :=
 x -\tau(z)$, the previous relation yields that
 $\dot {\tilde x}(t) = F \tilde x(t)$
 for all $t \geq 0$ and for all initial states $x_0 \in \Real^m$ and $z_0 \in
 {\ca}$. This, the fact that $F$ is Hurwitz and that $\ca$ is
 locally asymptotically (exponentially) stable for (\ref{zdbold}) yield immediately
 the desired result. $\triangleleft$

 \subsection{Proof of Proposition \ref{propositionvarrho}} \label{proofpropositionvarrho}
 The result will be proved by taking the ``complex" pair
 \beeq{\label{Fcompl}
 F= \mbox{diag}(\lambda_1, \ldots, \lambda_{r+1}) \qquad G=( g, \ldots,\, g )\tr
 }
 in which $\lambda_i \in \Compl_\ell=\{\lambda \in \Compl \; : \;
 \mbox{Re}\lambda<-\ell\}$, $i=1,\ldots, r+1$, $\ell>0$, and $g\neq 0$. Once proved the
 result for the $(r+1)$-dimensional pair in (\ref{Fcompl}), the claim of the
 proposition follows by taking any $(2r + 2)$-dimensional ``real" representation of
 (\ref{Fcompl}).

 By bearing in mind the definition of the map $\tau$ in (\ref{taumapdef})
 note that, as a consequence of the choice of $F$ and $G$ in (\ref{Fcompl}),
 it turns out that
  \beeq{ \label{taudiag}
 \tau(z) = \left (\ba{cccc}
 \tau_{\lambda_1}(z) &
 \tau_{\lambda_2}(z) &
 \cdots &
 \tau_{\lambda_{r+1}}(z)
 \ea \right )\tr
 \qquad
 \tau_{\lambda_i}(z) = \int_{-\infty}^0
 e^{- \lambda_i s} \,g\, {{q}_0}(\hat z(s, z)) ds\,.
 }
 We will prove next that there exists an $\ell>0$ such that
 by arbitrarily choosing $\lambda_i$, $i=1,\ldots, r+1$, in
 $\Compl_\ell \setminus S$, where $S$ is a set of zero Lebesgue measure,
 then the map $\tau$ is such that
 \beeq{\label{tauqcond}
 \tau(z_1) = \tau(z_2) \qquad \Rightarrow \qquad {{q}_0}(z_1) =
 {{q}_0}(z_2)
 \qquad \forall  \, z_1, z_2 \in \Real^n\,.
 }

 More precisely we will prove that, having defined
 \[
 \Upsilon = \{ (z_1, z_2) \in \Real^n \times \Real^n\; : \;
 {{q}_0}(z_1)\neq
 {{q}_0}(z_2) \}\,,
\]
 the set
 \beeq{\label{SetZeromeas}
 S = \{(\lambda_1, \ldots, \lambda_{r+1}) \in \Compl_\ell^{r+1}\quad : \quad
 \exists \; (z_1,z_2) \in \Upsilon \; : \; \tau_{\lambda_i}(z_1)=
 \tau_{\lambda_i}(z_2)
 \quad \forall \; i =1, \ldots, r+1 \}
 }
 has zero Lebesgue measure in $\Compl^{r+1}$ for a proper choice of $\ell$.
 To this end the following Theorem,
 proved in a more general setting
  in \cite{AnPr} (see also (\cite{Coron})),
 plays a crucial role.

 \begin{theorem} \label{TheoremCoron}
 Let $\Omega$ and $\Upsilon$ be open subsets of $\Compl$ and $\Real^{2n}$
 respectively. Let $(\varpi,\lambda )\in\Upsilon \times \Omega  \mapsto
\delta_\tau(\varpi,\lambda )\in\Compl$
 be a function which is
 holomorphic in $\lambda$ for each $\varpi \in\Upsilon$ and $C^1$ for each
 $\lambda \in \Omega$. If, for each pair $(\varpi,\lambda) \in \Upsilon \times \Omega$
 for which $\delta_\tau(\varpi,\lambda)$ is zero there exists an integer $k$ satisfying
 \beeq{ \label{CoronCondition}
 \ba{l}
 \dst {\partial^i \delta_\tau \over \partial \lambda^i} (\varpi,\lambda) =0 \qquad
 \mbox{for all } i \in \{0, \ldots, k-1 \}\\[3mm]
 \dst {\partial^k \delta_\tau \over \partial \lambda^k} (\varpi,\lambda) \neq 0
 \ea
 }
 then the set
 \[
   S = \bigcup_{\varpi\in \Upsilon}\{(\lambda_1, \ldots, \lambda_{r+1}) \in {\Omega}^{r+1}\quad : \quad
\delta_\tau(\varpi,\lambda_1)=
 \ldots =\delta_\tau(\varpi,\lambda_{r+1})=0\}
 \]
 has zero Lebesgue measure in $\Compl^{r+1}$.
 \end{theorem}

 To apply this theorem to our context we first observe the following.
 \begin{proposition} \label{Proptauc1}
  There exists an $\ell>0$ such that for all $\lambda_i \in \Compl_\ell $, $i=1,
  \ldots, r+1$, the map $\tau(\cdot)$ in (\ref{taudiag}) is $C^1$.
 \end{proposition}

\begin{proof}
The map $\tau(\cdot)$ in (\ref{taudiag}) is $C^1$ if the functions
$e^{-\lambda_i s} \, g
\, {\partial {{q}_0}(\hat z(s, z))/
\partial z}$, $i = 1,\ldots, r+1$, are integrable on $s \in (-\infty, 0]$ for all $z \in
\Real^n$ (see \cite{Flemm}).  Consider the expansion
\[
{\partial {{q}_0}(\hat z(s, z))\over
\partial z}=\Bigl[{\partial q_0 \over \partial z}\Bigr]_{[z=\hat z(s,
z)]}{\partial \hat z(s, z) \over \partial z}\,.\] By definition,
there is a number $M$ such that $|\hat z(s, z)|\le M$ for all
$s\le 0$ and all $z\in \Real^n$. This, along with the fact that
$q_0(z)$ is $C^1$, shows that the first factor is bounded on
$(-\infty,0]\times \Real^n$. As for the second factor, bearing in
mind the notation introduced in Section \ref{proofpropositiontau},
observe that
\[
{d \over ds}{\partial \hat z(s, z) \over \partial z}
=\Bigl[{\partial a_0(z) {f}_0(z) \over
\partial z} \Bigr]_{[z=\hat
z(s, z)]}{\partial
\hat z(s, z) \over \partial z}
\]
Letting
 \[
\bar f = \max_{z \in \Real^n}{\partial a_0(z) {f}_0(z) \over
\partial z}
 \]
we obtain
\[
|{\partial
\hat z(s, z) \over \partial z}| \leq e^{\bar f |s|}\]
 for all $s$ and
 for all $z \in \Real^n$. From this, the result immediately follows with $\ell=\bar f$.
\end{proof} $\triangleleft$ {\em (End proof Proposition \ref{Proptauc1})}

\vspace{5mm}

 Now set $\varpi:= (z_1, z_2)$ and

 \[
 \delta_\tau(\varpi,\lambda)= \int_{-\infty}^0 e^{-\lambda s} \, g\, \left [
 {{q}_0}(\hat z(s,z_1)) - {{q}_0}(\hat z(s,z_2))
 \right ]ds = \tau_\lambda(z_1) - \tau_\lambda(z_2)\,.
 \]
 This function is $C^1$ in $\varpi \in \Real^n \times \Real^n$ and it is holomorphic
 in $\lambda \in \Compl_\ell$ for every $\varpi\in \Real^n \times \Real^n$
 (see \cite{Rudin}, Chap. 19, p. 367). Moreover, as
 \[
 \int_{-\infty}^0 e^{-a s} \, \left | g \,{{q}_0}(\hat z(s,z_1)) -
 g \,{{q}_0}(\hat z(s,z_2)) \, \right |^2 ds < + \infty
 \]
 for all $\varpi \in \Real^n \times \Real^n$ and for all $a<0$,
 the Plancherel Theorem can be invoked to obtain
 \beeq{\label{Plancherel}
 {1 \over 2 \pi} \int_{-\infty}^\infty |\delta_\tau(\varpi, a+ i s)|^2 ds
 = \int_{-\infty}^0 e^{- 2 a s} \left | g \,{{q}_0}(\hat z(s,z_1)) -
 g \,{{q}_0}(\hat z(s,z_2)) \, \right |^2 ds
 }
 for all $a<0$ and for all $\varpi \in \Real^n \times \Real^n$.

   Now note that, for $\varpi=(z_1,z_2) \in \Upsilon$, we have $q_0(z_1) \neq
   q_0(z_2)$ and by continuity of flow with respect to time, there exists a time $t_1<0$ such that
   \[
   \left | g \,{{q}_0}(\hat z(s,z_1)) -
 g \,{{q}_0}(\hat z(s,z_2)) \right |>0 \qquad \mbox{for all } s \in (t_1, 0]
   \]
   which, combined with (\ref{Plancherel}), yields
   \[
    \int_{-\infty}^\infty |\delta_\tau(\varpi, a+ i s)|^2 ds >0\,.
   \]
   This implies that, for each $\varpi \in \Upsilon$, the function $\lambda \mapsto
   \delta_\tau(\varpi,\lambda)$ is not identically zero on $\Compl_\ell$. Since it is
   holomorphic, it turns out that for each $(\varpi,\lambda) \in \Upsilon \times
   \Compl_\ell$ it is possible to find an integer $k$ satisfying (\ref{CoronCondition}).
   Hence Theorem \ref{TheoremCoron} can be applied to obtain the desired
   result, namely that the set (\ref{SetZeromeas}) has zero Lebesgue measure.

   By this result we are guaranteed that by arbitrarily picking
   $r+1$ complex eigenvalues in $\Compl_\ell \setminus S$
   (with $\ell$ dictated by Proposition
   \ref{Proptauc1}) of  $F$ defined in
   (\ref{Fcompl}), condition (\ref{tauqcond}) is satisfied. From this it is
   easy to show that there exists a class-${\cal K}$ function satisfying
   (\ref{partialinjection}). For, define
   \[
   \varphi(s) = \sup_{\small \ba{l}
   |\tau(z_1) - \tau(z_2)| \leq s\\
   z_1,z_2 \in {\ca}\ea}
   |{{q}_0}(z_1) -  {{q}_0}(z_2)| \,.
   \]
  This function is increasing and, as a consequence of (\ref{tauqcond}),
  $\varphi(0)=0$. Moreover it is possible to prove that $\varphi(s)$ is
  continuous at $s=0$. For, suppose that it is not, namely, as
  $\varphi(\cdot)$ is increasing and $\varphi(0)=0$,
  suppose that there exists a $\varphi^\star>0$ such that
  $\lim_{s \to 0^+} \varphi(s) = \varphi^\star$. This implies that
  there exist sequences $\{z_{1n}\}$, $\{z_{2n}\}$ in $\ca$, such that
  $|{{q}_0}(z_{1n}) -  {{q}_0}(z_{2n})| \geq \varphi^\star/2$ and
  $|\tau({z_{1n}}) - \tau(z_{2n})| <1/n$ for any $n \in \Integer$.
  But, as $\ca$ is bounded, there are subsequences of $\{z_{1n}\}$,
  $\{z_{2n}\}$ which, for $n \to \infty$, converge to $z_{1}^\star$,
  $z_{2}^\star$ respectively. As $\tau(\cdot)$ and ${{q}_0}(\cdot)$ are
  continuous $\tau(z_{1}^\star) - \tau(z_{2}^\star)=0$ and
  $|{{q}_0}(z_{1}^\star) -  {{q}_0}(z_{2}^\star)| \geq \varphi^\star/2$
  which contradict (\ref{tauqcond}). Hence $\varphi(s)$ is continuous at
  $s=0$. With this result at hand, define the candidate class-$\cal K$ function
  \[
   \varrho(s) = {1 \over s} \int_s^{2s} \varphi(\sigma) d \sigma\;+\;
            s
  \]
 which satisfies
  \beeq{ \label{varrhoineq}
   \varphi(s) \leq \varrho(s) \,.
   }
  By construction this function is continuous for all $s >0$ and, as $\varphi(s)$
  is continuous at $s=0$ and by (\ref{varrhoineq}), it is also continuous at
  $s=0$. Moreover, by (\ref{varrhoineq}) and by definition of $\varphi(\cdot)$,
  (\ref{partialinjection}) is also satisfied.

 \subsection{Proof of Proposition \ref{Propositionfinal}}
 \label{proofpropositionfinal}

  By the result of the previous proposition we know that
 \[
 \tau(z_1) = \tau(z_2) \qquad \Rightarrow \qquad {{q}_0}(z_1) = {{q}_0}(z_2)
 \qquad \forall  \, z_1, z_2 \in \ca\,.
 \]
For any $x\in \tau({\cal A})$, let $[x]=\{z\in {\cal A}:
\tau(z)=x\}$. The previous property shows that the map
$q_0(\cdot)$ is constant on $[x]$. As a consequence, there is a
well defined function $\gamma_0: \tau({\cal A}) \to
\Real$ such that
\[
\gamma_0(\tau(z)) = - q_0(z), \qquad \forall z\in {\cal A}\,.
\]
In fact, the value $\gamma_0(x)$ at any  $x \in \tau({\cal A})$ is
simply defined by taking any $z\in [x]$ and setting
$\gamma_0(x):=-q_0(z)$.
  Moreover, by (\ref{partialinjection}), the map in
  question is also continuous. Now note that, as $\ca$ is compact and
  $\tau(\cdot)$ and ${{q}_0}(\cdot)$
  are continuous maps, $\tau({\ca}) \subset \Real^m$ and ${{q}_0}({\ca})
  \subset \Real$ are compact sets. From this the Tietze's extension
  Theorem (see, for instance, Theorem VII.5.1 in \cite{Dugundji}) can be invoked to
  claim the existence of a continuous map
 $\gamma: \Real^m \to \Real$  which agrees with $\gamma_0$ on $\tau({\cal
  A})$. This implies that ${{q}_0}(z) + \gamma \circ \tau(z)=0$
for all   $z \in \ca$ and proves the first claim of the
proposition.

 Furthermore if $\varrho(\cdot)$ is linearly bounded at the origin, by compactness arguments
 it is possible to claim the existence of a positive $\bar \varrho$ such that
 $|{{q}_0}(z_1) - {{q}_0}(z_2)| \leq \bar \varrho |  \tau(z_1) -
 \tau(z_2)|$. It follows that $\gamma_0$ is a Lipschitz function on $\tau(\ca)$.
 From this the Kirszbraun Theorem (see, for instance,
 Theorem 2.10.43 in \cite{Federer}) yields the existence of a Lipschitz map
 $\gamma:\Real^m \to \Real$, with Lipschitz constant $\varrho$, which agrees
 with $\gamma_0$ on $\tau({\ca})$. This completes the proof of the
 proposition.

\end{document}